\documentclass[global]{svjour}
\usepackage{amssymb}
\usepackage{amsfonts}
\usepackage{amsmath}

\input{tcilatex}

\begin{document}

\title{A Fixed Point Theorem and Equilibria of Abstract Economies with
w-Upper Semicontinuous Set-Valued Maps}
\author{Carlos Herv\'{e}s-Beloso$^{1}$ and Monica Patriche$^{2}$}
\institute{1 University of Vigo 
%TCIMACRO{\TeXButton{email}{\email{cherves@uvigo.es}}}%
%BeginExpansion
\email{cherves@uvigo.es}%
%EndExpansion
\\
2 University of Bucharest 
%TCIMACRO{\TeXButton{email}{\email{monica.patriche@yahoo.com}}}%
%BeginExpansion
\email{monica.patriche@yahoo.com}%
%EndExpansion
}
\mail{1 RGEA, Facultad de Econ\'{o}micas, Universidad de Vigo, Campus
Universitario, E-36310 Vigo, Spain \\
2 University of Bucharest, Faculty of Mathematics and Computer Science, 14
Academiei Street, 010014 Bucharest, Romania}
\maketitle

\begin{abstract}
{\small We introduce the notion of w-upper semicontinuous set valued maps }
\end{abstract}


{\small \ }

\bigskip

\noindent \textbf{1.} \textbf{Introduction\medskip }

\noindent The pioneer work of Nash [11] first proved a theorem of
equilibrium existence for games where the player's payoffs are represented
by continuous quasi-concave utilities. Arrow and Debreu used the work by
Nash to prove the existence of equilibrium in a generalized N-person game or
on abstract economy [7] which implies the Walrasian equilibrium existence
[2]. These ideas were extended by various authors in several ways. In [16],
Shafer and Sonnenschein proved the existence of equilibrium of an economy
with finite dimensional commodity space and irreflexive preferences
represented as set valued maps with open graph. Yannelis and Prahbakar [22]
developed new techniques based on selection theorems and fixed-point
theorems. Their main result concerns the existence of equilibrium when the
constraint and preference set valued maps have open lower sections. They
work within different frameworks (countable infinite number of agents,
infinite dimensional strategy spaces).

Borglin and Keiding [3] used new concepts of K.F.-set valued maps and
KF-majorized set valued maps for their existence results . The concept of
KF-majorized set valued maps was extended by Yannelis and Prabhakar [22] to
L-majorized set valued maps. In [23], Yuan proposed a more general model of
abstract economy than the one introduced by Borglin and Keing in [3], in the
sense that the constraint mapping was split into two parts $A$ and $B.$ This
is due to the "small" constraint set valued map $A$ which could not have
enough fixed points even though the "big" constraint set valued map $B$
could.

Most existence theorems of equilibrium deal with preference set valued maps
which have lower open sections or are majorized by set valued maps with
lower open sections. In the last few years, some existence results were
obtained for lower semicontinuous and upper semicontinuous set valued maps.
Some recent results concerning upper semicontinuous set valued maps and
fixed points can be found in [1], [4], [18], [19], [20], [24]. New results
on equilibrium existence in games are given in [10], [12], [13], [17].

In this paper, we define two types of set valued maps: w-upper
semicontinuous set valued maps and set valued maps that have e-USS-property.
We prove a fixed point theorem for \textit{w-}upper semicontinuous set
valued maps. This result is a Wu like result [20] and generalizes the
Himmelberg's fixed point theorem in [9]. We use this theorem for proving our
first theorem of equilibrium existence for abstract economies having w-upper
semicontinuous constraint and preference set valued maps. On the other hand,
we use a technique of approximation to prove an equilibrium existence
theorem for set valued maps having e-USS-property.

The paper is organized in the following way: Section 2 contains
preliminaries and notations. The fixed point theorem is presented in Section
3 and the equilibrium theorems are stated in Section 4.\bigskip

\noindent \textbf{2.} \textbf{Preliminaries and Notation}\medskip

\noindent Throughout this paper, we shall use the following notations and
definitions:

Let $A$ be a subset of a topological space $X$. $2^{A}$ denotes the family
of all subsets of $A$. cl$A$ denotes the closure of $A$ in $X$. If $A$ is a
subset of a vector space, co$A$ denotes the convex hull of $A$. If $F$, $G:$ 
$X\rightrightarrows Y$ are set valued maps, then conv $G$, cl $G$, $G\cap F$ 
$:$ $X\rightrightarrows Y$ are set valued maps defined by $($conv $G)(x)=$%
conv $G(x)$, $($cl $G)(x)=$cl $G(x)$ and $(G\cap F)(x)=G(x)\cap F(x)$ for
each $x\in X$, respectively. The graph of $T:X\rightrightarrows Y$ is the
set Gr $(T)=\{(x,y)\in X\times Y\mid y\in T(x)\}.$

The set valued map $\overline{T}$ is defined by $\overline{T}(x):=\{y\in
Y:(x,y)\in $cl$_{X\times Y}$ Gr $T\}$ (the set cl$_{X\times Y}$ Gr $(T)$ is
called the adherence of the graph of $T$)$.$ It is easy to see that cl $%
T(x)\subset \overline{T}(x)$ for each $x\in X.\medskip $

Let $X$, $Y$ be topological spaces and $T:X\rightrightarrows Y$ be a set
valued map. $T$ is said to be \textit{upper semicontinuous} if for each $%
x\in X$ and each open set $V$ in $Y$ with $T(x)\subset V$, there exists an
open neighborhood $U$ of $x$ in $X$ such that $T(y)\subset V$ for each $y\in
U$. $T$ is said to be \textit{almost upper semicontinuous} if for each $x\in
X$ and each open set $V$ in $Y$ with $T(x)\subset V$, there exists an open
neighborhood $U$ of $x$ in $X$ such that $T(y)\subset $cl $V$ for each $y\in
U$.\medskip

\noindent \textbf{Lemma 2.1}(Lemma 3.2, pag. 94 in [25]) \textit{Let }$X$%
\textit{\ be a topological space, }$Y$\textit{\ be a topological linear
space, and let }$S:X\rightrightarrows Y$\textit{\ be an upper semicontinuous
set valued map with compact values. Assume that the set }$C\subset Y$\textit{%
\ is closed and }$K\subset Y$\textit{\ is compact. Then }$%
T:X\rightrightarrows Y$\textit{\ defined by }$T(x)=(S(x)+C)\cap K$\textit{\
for all }$x\in X$\textit{\ is upper semicontinuous.}\medskip

Lemma 2.2 is a version of Lemma 1.1 in [21] ( for $D=Y,$ we obtain Lemma 1.1
in [21]). For the reader's convenience, we include its proof below.

\noindent \textbf{Lemma 2.2} \textit{Let }$X$\textit{\ be a topological
space, }$Y$\textit{\ be a nonempty subset of a locally convex topological
vector space }$E$\textit{\ and }$T:X\rightrightarrows Y$\textit{\ be a set
valued map}$.$\textit{\ Let \ss\ be a basis of neighbourhoods of }$0$\textit{%
\ in }$E$\textit{\ consisting of open absolutely convex symmetric sets. Let }%
$D$\textit{\ be a compact subset of }$Y$\textit{. If for each }$V\in $%
\textit{\ss , the set valued map }$T^{V}:X\rightrightarrows Y$\textit{\ is
defined by }$T^{V}(x)=(T(x)+V)\cap D$\textit{\ for each }$x\in X,$\textit{\
then }$\cap _{V\in \text{\ss }}\overline{T^{V}}(x)\subseteq \overline{T}(x)$%
\textit{\ for every }$x\in X.\medskip $

\noindent \textit{Proof} Let $x$ and $y$ be such that $y\in \cap _{V\in 
\text{\ss }}\overline{T^{V}}(x)$ and suppose, by way of contradiction, that $%
y\notin \overline{T}(x).$ This means that $(x,y)\notin $cl Gr $T,$ so that
there exists an open neighborhood $U$ of $x$ and $V\in $\ss\ such that:

$(U\times (y+V))\cap $Gr $T=\emptyset .\ \ \ \ \ \ \ \ \ \ \ \ \ \ \ \ \ \ \
\ \ \ \ \ \ \ \ \ \ \ \ \ (1)$

Choose $W\in $\ss\ such that $W-W\subseteq V$ (e.g. $W=\frac{1}{2}V)$. Since 
$y\in T^{W}(x)$, then $(x,y)\in $cl Gr $T^{W},$ so that

\begin{equation*}
(U\times (y+W))\cap \text{Gr }T^{W}\neq \emptyset .\ \ \ \ \ \ \ \ \ \ \ \ \
\ \ \ \ \ \ \ \ \ \ \ \ \ \ \ \ \ \ \ 
\end{equation*}

There are some $x^{\prime }\in U$ and $w^{\prime }\in W$ such that $%
(x^{\prime },y+w^{\prime })\in $Gr $T^{W},$ i.e. $y+w^{\prime }\in
T^{W}(x^{\prime }).$ Then, $y+w^{\prime }\in D$ and $y+w^{\prime }=y^{\prime
}+w^{^{\prime \prime }}$ for some $y^{\prime }\in T(x^{\prime })$ and $%
w^{^{\prime \prime }}\in W.$ Hence, $y^{\prime }=y+(w^{\prime }-w^{^{\prime
\prime }})\in y+(W-W)\subseteq y+V,$ so that $T(x^{\prime })\cap (y+V)\neq
\emptyset .$ Since $x^{\prime }\in U,$ this means that $(U\times (y+V))\cap $%
Gr $T\neq \emptyset ,$ contradicting (1). $\ \ \ \ \ \ \ \ \ \ \square $%
\medskip

We introduce the following definitions.

Let $X$ be a topological space, $Y$ be a nonempty subset of a topological
vector space $E$ and $D$ be a subset of $Y$.

\noindent \textbf{Definition 2.1} The set valued map $T:X\rightrightarrows Y$
is said to be \textit{w-upper semicontinuous} (weakly upper semicontinuous) 
\textit{with respect to} \textit{the set} $D$ if there exists a basis \ss\ %
of open symmetric neighborhoods of $0$ in $E$ such that, for each $V\in $\ss %
, the set valued map $T^{V}$ is upper semicontinuous.

\noindent \textbf{Definition 2.2 }The set valued map $T:X\rightrightarrows Y$
is said to be \textit{almost} \textit{w-upper semicontinuous} (almost weakly
upper semicontinuous) \textit{with respect to} \textit{the set} $D$ if there
exists a basis \ss\ of open symmetric neighborhoods of $0$ in $E$ such that,
for each $V\in $\ss , the set valued map $\overline{T^{V}}$ is upper
semicontinuous.\medskip

\noindent \textit{Example}\textbf{\ }\textit{2.1} Let $T_{1}:(0,2)%
\rightrightarrows (0,2)$ be defined by $T_{1}(x):=\left\{ 
\begin{array}{c}
(0,1)\text{ if }x\in (0,1]; \\ 
\lbrack 1,2)\text{ if x}\in (1,2).%
\end{array}%
\right. $

$T_{1}$ and $T_{1}\cap \{1\}:=\left\{ 
\begin{array}{c}
\phi \text{ \ if \ }x\in (0,1]; \\ 
\{1\}\text{ if x}\in (1,2)%
\end{array}%
\right. $ are not upper semicontinuous on $(0,2).$

Let $D:=\{1\}$ and let $V:=(-\varepsilon ,\varepsilon ),$ $\varepsilon >0,$
be an open symmetric neighbourhood of $0$ in $IR.$ Then, it results that

for $\varepsilon >0,$ $T_{1}(x)+(-\varepsilon ,\varepsilon ):=\left\{ 
\begin{array}{c}
(-\varepsilon ,1+\varepsilon )\text{ \ \ if \ \ \ }x\in (0,1]; \\ 
(1-\varepsilon ,2+\varepsilon )\text{\ if \ \ \ \ \ }x\in (1,2);%
\end{array}%
\right. $

$T_{1}^{V}(x):=(T_{1}(x)+(-\varepsilon ,\varepsilon ))\cap \{1\}=\{1\}$ for
any $x\in (0,2).$

$\overline{T_{1}^{V}}(x)=\{1\}$ for $x\in (0,2).$

For each $V=(-\varepsilon ,\varepsilon )$ with $\varepsilon >0,$ the set
valued maps $T_{1}^{V}$ and $\overline{T_{1}^{V}}$ are upper semicontinuous
and $\overline{T_{1}^{V}}$ has nonempty values. We conclude that $T_{1}$ is
w-upper semicontinuous with respect to $D$ and it is also almost w-upper
semicontinuous with respect to $D.$

We also define the dual w-upper semicontinuity with respect to a compact set.

\noindent \textbf{Definition 2.3 }Let $T_{1},T_{2}:X\rightrightarrows Y$ be
set valued maps. The pair $(T_{1},T_{2})$ is said to be \textit{dual almost
w-upper semicontinuous} (dual almost weakly upper semicontinuous) \textit{%
with respect to} the set $D$ if there exists a basis \ss\ of open symmetric
neighborhoods of $0$ in $E$ such that, for each $V\in $\ss , the set valued
map $\overline{T_{(1,2)}^{V}}:X\rightrightarrows D$ is lower semicontinuous,
where $T_{(1,2)}^{V}:X\rightrightarrows D$ is defined by $%
T_{(1,2)}^{V}(x):=(T_{1}(x)+V)\cap T_{2}(x)\cap D$ for each $x\in X$.\medskip

\noindent \textit{Example}\textbf{\ }\textit{2.2} Let $\ D:=[1,2],$ $%
T_{1}:(0,2)\rightrightarrows \lbrack 1,4]$ be the set valued map defined by

$T_{1}(x):=\left\{ 
\begin{array}{c}
\lbrack 2-x,2],\text{ if }x\in (0,1); \\ 
\{4\}\text{ \ \ \ \ \ \ if \ \ \ \ \ \ \ }x=1; \\ 
\lbrack 1,2]\text{ \ \ \ if \ \ \ }x\in (1,2).%
\end{array}%
\right. $

and $T_{2}:(0,2)\rightrightarrows \lbrack 2,3]$ be the set valued map
defined by

$T_{2}(x):=\left\{ 
\begin{array}{c}
\lbrack 2,3],\text{ if }x\in (0,1]; \\ 
\{2\}\text{ \ if \ \ }x\in (1,2);%
\end{array}%
\right. .$

The set valued map $T_{1}$ is not upper semicontinuous on $(0,2)$.

For $\varepsilon \in (0,2],$ $(T_{1}(x)+(-\varepsilon ,\varepsilon ))\cap
D\cap T_{2}(x)=\left\{ 
\begin{array}{c}
\{2\}\text{ if }x\in (0,1)\cup (1,2); \\ 
\phi \text{ \ \ \ \ \ \ \ if \ \ \ \ \ \ \ \ \ }x=1.%
\end{array}%
\right. $

For $\varepsilon \in (2,\infty ),$ $(T_{1}(x)+(-\varepsilon ,\varepsilon
))\cap D\cap T_{2}(x)=\{2\}$ for each $x\in (0,2).$

Then, we have that for each $\varepsilon >0,$ $\overline{T_{(1,2)}^{V}}%
(x)=\{2\}$ for each $x\in \lbrack 0,2]$ and the set valued map $\overline{%
T_{(1,2)}^{V}}$ is upper semicontinuous and has nonempty values.

We conclude that the pair $(T_{1},T_{2})$ is dual almost w-upper
semicontinuous with respect to $D.$\medskip

\noindent \textbf{3.} \textbf{A New Fixed Point Theorem\medskip }

\noindent We obtain the following fixed point theorem which generalizes
Himmelberg's fixed point theorem in [9]:

\noindent \textbf{Theorem 3.1} \textit{Let }$I$\textit{\ be an index set.
For each }$i\in I,$\textit{\ let }$X_{i}$\textit{\ be a nonempty convex
subset of a Hausdorff locally convex topological vector space }$E_{i}$%
\textit{, }$D_{i}$\textit{\ be a nonempty compact convex subset of }$X_{i}$%
\textit{\ and }$S_{i},T_{i}:X:=\tprod\limits_{i\in I}X_{i}\rightrightarrows
X_{i}$\textit{\ be two set valued maps with the following conditions:}

\textit{1) for each }$x\in X$\textit{, }$\overline{S}_{i}(x)\subseteq
T_{i}(x)$\textit{; }

\textit{2)\ }$S_{i}$\textit{\ is almost w-upper semicontinuous with respect
to }$D_{i}$\textit{\ and }$\overline{S_{i}^{V_{i}}}$\textit{\ is convex
nonempty valued for each absolutely convex symmetric neighborhood }$V_{i}$%
\textit{\ of }$0$\textit{\ in }$E_{i}$\textit{.}

\textit{Then there exists }$x^{\ast }\in D:=\tprod\limits_{i\in I}D_{i}$%
\textit{\ such that }$x_{i}^{\ast }\in T_{i}(x^{\ast })$\textit{\ for each }$%
i\in I.\medskip $

\noindent \textit{Proof }Since $D_{i}$ is compact, $D:=\tprod\limits_{i\in
I}D_{i}$ is also compact in $X.$ For each $i\in I,$ let \ss $_{i}$\ be a
basis of open absolutely convex symmetric neighborhoods of zero in $E_{i}$
and let \ss =$\tprod\limits_{i\in I}$\ss $_{i}.$ For each system of
neighborhoods $V=(V_{i})_{i\in I}\in \tprod\limits_{i\in I}$\ss $_{i},$
let's define the set valued maps $S_{i}^{V_{i}}:X\rightrightarrows D_{i},$
by $S_{i}^{V_{i}}(x)=(S_{i}(x)+V_{i})\cap D_{i}$, $x\in X,$ $i\in I.$ By
assumption 2) each $\overline{S_{i}^{V_{i}}}$ is u.s.c with nonempty closed
convex values. Let's define $S^{V}:X\rightrightarrows D$ by $%
S^{V}(x)=\tprod\limits_{i\in I}\overline{S_{i}^{V_{i}}}(x)$ for each $x\in
D. $ The set valued map $S^{V}$ is upper semicontinuous with closed convex
values. Therefore, according to Himmelberg's fixed point theorem [9], there
exists $x_{V}^{\ast }=\tprod\limits_{i\in I}(x_{V}^{\ast })_{i}\in D$ such
that $x_{V}^{\ast }\in S^{V}(x_{V}^{\ast }).$ It follows that $(x_{V}^{\ast
})_{i}\in \overline{S_{i}^{V_{i}}}(x_{V}^{\ast })$ for each $i\in I.$

For each $V=(V_{i})_{i\in I}\in $\ss $,$ let's define $Q_{V}=\cap _{i\in
I}\{x\in D:$ $x_{i}\in \overline{S_{i}^{V_{i}}}(x)\}.$

$Q_{V}$ is nonempty since $x_{V}^{\ast }\in Q_{V},$ then $Q_{V}$ is nonempty
and closed.

We prove that the family $\{Q_{V}:V\in \text{\ss }\}$ has the finite
intersection property.

Let $\{V^{(1)},V^{(2)},...,V^{(n)}\}$ be any finite set of $\text{\ss }$ and
let $V^{(k)}=\underset{i\in I}{\tprod }V_{i}^{(k)}$, $k=1,...,n.$ For each $%
i\in I$, let $V_{i}=\underset{k=1}{\overset{n}{\cap }}V_{i}^{(k)}$, then $%
V_{i}\in \text{\ss }_{i};$ thus $V=\underset{i\in I}{\tprod }V_{i}\in 
\underset{i\in I}{\tprod }\text{\ss }_{i}.$ Clearly $Q_{V}\subseteq \underset%
{k=1}{\overset{n}{\cap }}Q_{V^{(k)}}$ so that $\underset{k=1}{\overset{n}{%
\cap }}Q_{V^{(k)}}\neq \emptyset .$

Since $D$ is compact and the family $\{Q_{V}:V\in \text{\ss }\}$ has the
finite intersection property, we have that $\cap \{Q_{V}:V\in \text{\ss }%
\}\neq \emptyset .$ Take any $x^{\ast }\in \cap \{Q_{V}:V\in $\ss $\},$ then
for each $V_{i}\in \text{\ss }_{i},$ $x_{i}^{\ast }\in \overline{%
S_{i}^{V_{i}}}(x^{\ast })$. According to Lemma 2.2,\emph{\ }we have that%
\emph{\ } $x_{i}^{\ast }\in \overline{S_{i}}(x^{\ast }),$ for each $i\in I,$
therefore $x_{i}^{\ast }\in T(x^{\ast }).$ $\ \ \ \ \ \ \ \ \ \ \ \ \ \ \ \
\ \ \ \ \ \ \ \ \ \ \ \ \ \ \ \ \ \ \ \ \ \ \ \ \ \ \ \ \ \ \ \ \ \ \ \ \ \
\ \ \ \ \ \ \ \ \ \ \ \ \ \ \ \ \ \ \ \square $\medskip

If $\left\vert I\right\vert =1$ we get the result below.

\noindent \textbf{Corollary 3.1} \textit{Let }$X$\textit{\ be a nonempty
subset of a Hausdorff locally convex topological vector space }$F,$\textit{\ 
}$D$\textit{\ be a nonempty compact convex subset of }$X$\textit{\ and }$%
S,T:X\rightrightarrows X$\textit{\ be two set valued maps with the following
conditions:}

\textit{1) for each }$x\in X,$\textit{\ }$\overline{S}(x)\subseteq T(x)$%
\textit{\ and }$S(x)\neq \emptyset ,$

\textit{2)\ }$S$\textit{\ is almost w-upper semicontinuous with respect to }$%
D$\textit{\ and }$\overline{S^{V}}$\textit{\ is convex valued for each open
absolutely convex symmetric neighborhood }$V$\textit{\ of }$0$\textit{\ in }$%
E$\textit{.}

\textit{Then, there exists a point }$x^{\ast }\in D$\textit{\ such that }$%
x^{\ast }\in T(x^{\ast }).$\medskip

In the particular case that the set valued map $S=T$ the following result
stands.

\noindent \textbf{Corollary 3.2 }\textit{Let }$X$\textit{\ be a nonempty
subset of a Hausdorff locally convex topological vector space }$F,$\textit{\ 
}$D$\textit{\ be a nonempty compact convex subset of }$X$\textit{\ and }$%
T:X\rightrightarrows X$\textit{\ be an almost w- upper semicontinuous set
valued map with respect to }$D$\textit{\ and }$\overline{T^{V}}$\textit{\ is
convex valued for each open absolutely convex symmetric neighborhood }$V$%
\textit{\ of }$0$\textit{\ in }$E$\textit{. Then, there exists a point }$%
x^{\ast }\in D$\textit{\ such that }$x^{\ast }\in \overline{T}(x^{\ast
}).\medskip $

\noindent \textbf{4. Application in the Equilibrium Theory\medskip }

\noindent \textbf{4.1 The Model of an Abstract Economy}

We will consider further Yuan's model of an abstract economy (see [23]). Let 
$I$ be a nonempty set (the set of agents). For each $i\in I$, let $X_{i}$ be
a non-empty subset of a topological vector space representing the agent's $i$
set of actions and define $X:=\underset{i\in I}{\prod }X_{i}$; let $A_{i}$, $%
B_{i}:X\rightrightarrows X_{i}$ be the constraint set valued maps and $P_{i}$
the preference set valued map.

\noindent \textbf{Definition 4.1}\textit{\ }[23]\textit{\ }An \textit{%
abstract economy} $\Gamma =(X_{i},A_{i},P_{i},B_{i})_{i\in I}$ is defined as
a family of ordered quadruples $(X_{i},A_{i},P_{i},B_{i})$.\medskip

\noindent \textbf{Definition 4.2 [23] }An \textit{equilibrium} for $\Gamma $
is defined as a point $x^{\ast }\in X$ such that for each $i\in I$, $%
x_{i}^{\ast }\in \overline{B}_{i}(x^{\ast })$ and $A_{i}(x^{\ast })\cap
P_{i}(x^{\ast })=\emptyset $.\medskip

\noindent \textit{Remark}\textbf{\ }\textit{4.1} When, for each $i\in I$, $%
A_{i}(x)=B_{i}(x)$ for all $x\in X,$ the abstract economy model coincides
with the classical one introduced by Borglin and Keiding in [3]. If in
addition, $\overline{B}_{i}(x^{\ast })=$cl$_{X_{i}}$ $B_{i}(x^{\ast })$ for
each $x\in X,$ which is the case where $B_{i}$ has a closed graph in $%
X\times X_{i}$, the definition of equilibrium coincides with the one used by
Yannelis and Prabhakar in [22].\medskip

\noindent \textit{Remark}\textbf{\ }\textit{4.2} If the preference set
valued map $P_{i}$ is defined by using a utility function $u_{i}$, that is $%
P_{i}(x)=\{y\in X_{i}:u_{i}(y)>u_{i}(x_{i})\},$ the irreflexibility
condition $x_{i}\notin \overline{P_{i}}\left( x\right) ,$ which appears
among the hypothesis of the existence equilibrium theorems, may fail. A case
in which this condition is verified, is when $P_{i}$ is an order interval
preference. Order interval preferences are studied, for instance, in
Chateauneuf [5]. These preference relations $\prec $ (on $X$) are
representable, if two real valued functions $u$ and $v$ on $X$ exist and are
such that: $x\prec y$ $\Leftrightarrow $ $u(x)<v(y)$. If a representation of
the preference relation $\prec _{i}$ exist, we can define the preference set
valued map $P_{i}$ by $P_{i}(x)=\{y\in X_{i}:v_{i}(y)>u_{i}(x_{i})\}$ and
the condition $x_{i}\notin \overline{P_{i}}\left( x\right) $ can be
fulfilled.\medskip

\noindent \textbf{4.2} \textbf{Examples of Abstract Economies with Two
Constraint Set Valued Maps}

\noindent A first example of an abstract economy with two constraint set
valued maps is the one associated to the model proposed by Radner [14] and
his followers. This is a model of a pure exchange economy with two periods,
present and future, and uncertainty on the state of nature in the future.
There is a finite number $n$ of agents and a finite number $m$ of possible
future states of nature. Let $I=\{1,2,...,n\}$ be the set of agents and $%
\Omega =\{s_{1},s_{2},...,s_{m}\}$ the set of the future states of the
nature. Each agent has his own private, and tipically incomplete,
information about the future state of nature. For each agent $i$ in $I$, the
initial private information is a partition on $\Omega $ , induced by a
signal $\pi _{i}:\Omega \rightarrow Y_{i}$. In Radner [14], agents make
decisions today without knowing the future state of nature tomorrow. The
initial agent's information is kept fixed and their consumption plans need
to be made compatible with their information, in the sense that their
consumption must be the same in states that they do not distinguish. In a
different framework, Radner [15] consider the notion of rational expectation
equilibrium. In this model, agents are able to forecast the future
equilibrium price. Consequently, their initial information is updated with a
signal given by the future equilibrium prices $p$ and a more refined
partition of $\Omega $ is obtained as the joint of the initial information
and the information generated by $\widehat{\pi _{i}}(p):\Omega \rightarrow
\Delta $, defined by $\widehat{\pi _{i}}(p)(s)=p_{i}(s)$, where $\Delta $ be
the normalized set of prices. Here, we consider a model in which agents may
be able to learn from market signals. These market signal are summarized by
equilibrium prices that may not be fully revealing. Without loss of
generality, let denote the joint of the initial information and the
information generated by prices by $\widehat{\pi _{i}}$.

For agent $i$ in $I$, the consumption plan in the first period will be
denoted by $x_{0}^{i}\in IR_{+}^{l}$ and in the second period, for each
state $s_{j},$ $j=1,2,...,m$, it will be denoted by $x_{j}^{i}\in
IR_{+}^{l}. $ A bundle for agent $i$ is $\
x^{i}=(x_{0}^{i},x_{1}^{i},x_{2}^{i},...,x_{m}^{i}).$ Let $%
X_{i}=IR_{+}^{ml+1}.$ Each agent has a preference set valued map $%
Q_{i}^{\prime }:\tprod\limits_{i\in I}X_{i}\rightrightarrows X_{i}$ and an
initial endowment $e^{i}=(e_{0}^{i},e_{1}^{i},e_{2}^{i},...,e_{m}^{i})\in
X_{i}.$

\noindent \textbf{Definition 4.3 }A \textit{pure exchange economy with
assymmetric information} is the family $\mathcal{E=}(I,\Omega ,\widehat{\pi }%
_{i},Q_{i}^{\prime },e^{i})_{i\in I}.$

\noindent \textbf{Definition 4.4 }An \textit{allocation} for the economy $%
\mathcal{E}$ is $x=(x^{i})_{i\in I}.$ The allocation is called \textit{%
phisically feasible} if $\tsum\limits_{i\in I}x^{i}\leq \tsum\limits_{i\in
I}e^{i}$ and \textit{informationally feasible for each agent }$i$ if $%
\widehat{\pi }_{i}(p)(s)=\widehat{\pi }_{i}(p)(s^{\prime })$ implies $%
x_{s}^{i}=x_{s^{\prime }}^{i}.$

Let $p_{0}$ be the price in the first period, for the second period let $%
p_{j}$ be the price in the state $j,$ $j=1,2,...,m$ and let $%
p=(p_{0},p_{1},...,p_{m}).$ Let $\Delta $ be the normalized set of prices.
Without loss of generality, we assume that $p$ belongs to $\Delta $.

The budget set valued map of agent $i$ is $B_{i}:\Delta \rightrightarrows
IR_{+}^{lm+1},$ defined by $B_{i}(p)=\{x^{i}\in
IR_{+}^{lm+1}:px^{i}<pe^{i}\}.$

The information set valued map of agent $i$ is $I_{i}:\Delta
\rightrightarrows IR_{+}^{lm+1},$ defined by $I_{i}(p)=\{x^{i}\in
IR_{+}^{lm+1}:x_{s}^{i}=x_{s^{\prime }}^{i}$ if $\widehat{\pi }_{i}(p)(s)=%
\widehat{\pi }_{i}(p)(s^{\prime })\}.$

\noindent \textbf{Definition 4.5 }The pair $(x^{\ast },p^{\ast })\in
IR_{+}^{n(lm+1)}\times \Delta $ is an\textit{\ equilibrium} for the
asymmetrically informed economy $\mathcal{E}$ if

1) $\tsum\limits_{i\in I}(x^{\ast })^{i}\leq \tsum\limits_{i\in I}e^{i}$

and for each $i\in I,$

2) $(x^{\ast })^{i}\in \overline{I_{i}}(p^{\ast })\cap \overline{B_{i}}%
(p^{\ast });$

3) $y^{i}\in Q_{i}^{\prime }(x^{\ast })\cap I_{i}(p^{\ast })$ implies that $%
y^{i}\notin B_{i}(p^{\ast }).\medskip $

Let $X:=\tprod\limits_{i\in I}X_{i}\times \Delta ,$ where for $i\in I,$ $%
X_{i}=IR_{+}^{lm+1}$ is the consumption set of agent $i.$ Let's define the
following set valued maps:

-for each $i\in I,$ $Q_{i}:X\rightrightarrows X_{i}$ is the preference set
valued map defined by $Q_{i}(x,p)=Q_{i}^{\prime }(x)$ for each $(x,p)\in X;$

- $Q_{n+1}:X\rightrightarrows \Delta $ is the preference set valued map
defined by $Q_{n+1}(x,p):=\{q\in \Delta :q(\tsum\limits_{i\in
I}(x^{i}-e^{i}))>p(\tsum\limits_{i\in I}(x^{i}-e^{i}))\}$ for each $(x,p)\in
X;$

-for $i\in I,$ $A_{i}:X\rightrightarrows 2^{X_{i}}$ is defined by $%
A_{i}(x,p):=\{y^{i}\in IR_{+}^{lm+1}:py^{i}<pe^{i}\}$ for each $(x,p)\in X;$

-$A_{n+1}:X\rightrightarrows \Delta $ is defined by $A_{n+1}(x,p):=\Delta $
for each $(x,p)\in X;$

-for $i\in I,$ $I_{i}:X\rightrightarrows X_{i}$ is defined by $%
I_{i}(x,p):=\{y^{i}\in IR_{+}^{lm+1}:y_{s}^{i}=y_{s^{\prime }}^{i}$ if $%
\widehat{\pi }_{i}(p)(s)=\widehat{\pi }_{i}(p)(s^{\prime })\}$ for each $%
(x,p)\in X;$

- $I_{n+1}:X\rightrightarrows \Delta $ is defined by $I_{n+1}(x,p):=\Delta $
for each $(x,p)\in X;\medskip $

\noindent \textbf{Definition 4.6 }The abstract economy associated to the
model of the pure exchange economy with assymmetric information is $\Gamma
=(X_{i},A_{i},P_{i},B_{i})_{i\in \{1,2,...,n+1\}},$ where:

-for $i\in I,$ $X_{i}:=IR_{+}^{lm+1}$ is the consumption set of agent $i$
and let $X:=\tprod\limits_{i\in I}X_{i}\times \Delta ;$

-$P_{i}:X\rightrightarrows X_{i}$ $(i\in I)$ and $P_{n+1}:X\rightrightarrows
\Delta $ are the preference set valued maps defined by $%
P_{i}(x,p)=Q_{i}(x,p)\cap I_{i}(x,p)$ for each $(x,p)\in X$ and $i\in
\{1,2,...,n+1\};$

-$A_{i}:X\rightrightarrows X_{i}$ $(i\in I)$ and $A_{n+1}:X\rightrightarrows
\Delta $ are the constraint set valued maps defined above$;$

-$B_{i}:X\rightrightarrows X_{i}$ $(i\in I)$ and $B_{n+1}:X\rightrightarrows
\Delta $ are the constraint set valued maps defined by $%
B_{i}(x,p):=A_{i}(x,p)\cap I_{i}(x,p)$ for each $(x,p)\in X$ and $i\in
\{1,2,...,n+1\}.\medskip $

\noindent \textit{Remark}\textbf{\ }\textit{4.3} We note that $%
A_{i}(x,p)\cap P_{i}(x,p)\subseteq B_{i}(x,p)$ for each $(x,p)\in X$ and for
each $i\in \{1,2,...,n+1\}.$

\noindent \textbf{Proposition 4.1} \textit{An equilibrium for the associated
abstract economy }$\Gamma $\textit{\ is an equilibrium of the economy with
assymmetric information }$E$\textit{.}

\textit{Proof Let }$(x^{\ast },p^{\ast })$ be an equilibrium for $\Gamma .$

1) For each $i\in \{1,2,...,n\},$ we have that $(x^{\ast })^{i}\in \overline{%
B_{i}}(x^{\ast },p^{\ast })=\overline{(A_{i}\cap I_{i})}(x^{\ast },p^{\ast
}) $ and then, by definition of $A_{i}$ and $I_{i},$ $(x^{\ast })^{i}\in 
\overline{(I_{i}\cap B_{i})}(p^{\ast })$;

2) $p^{\ast }\in \overline{B_{n+1}}(x^{\ast },p^{\ast })=\Delta ;$

3) for each $i\in \{1,2,...,n\},$ we have that $A_{i}(x^{\ast },p^{\ast
})\cap P_{i}(x^{\ast },p^{\ast })=\phi $, which implies that if $y^{i}\in
P_{i}(x^{\ast },p^{\ast })=Q_{i}(x^{\ast },p^{\ast })\cap I_{i}(x^{\ast
},p^{\ast }),$ then $y^{i}\notin A_{i}(x^{\ast },p^{\ast });$ This means
that $y^{i}\in Q_{i}^{\prime }(x^{\ast })\cap I_{i}(p^{\ast })$ implies that 
$y^{i}\notin B_{i}(p^{\ast });$

4) we have that $A_{n+1}(x^{\ast },p^{\ast })\cap P_{n+1}(x^{\ast },p^{\ast
})=\phi $, which is equivalent with $\{q\in \Delta :q(\tsum\limits_{i\in
I}((x^{\ast })^{i}-e^{i}))>p^{\ast }(\tsum\limits_{i\in I}((x^{\ast
})^{i}-e^{i}))\}\cap \Delta =\phi .$ This fact implies that $%
q(\tsum\limits_{i\in I}((x^{\ast })^{i}-e^{i}))\leq p^{\ast
}(\tsum\limits_{i\in I}((x^{\ast })^{i}-e^{i}))\leq 0$ for all $q\in \Delta
. $ If we choose $q$ as a vector of the canonical basis of $IR^{ml+1},$ that
is $q_{j}=1$ and $q_{i}=0$ for $i\neq j,$ where $i,j\in \{1,2,...,ml+1\},$
we obtain that $\tsum\limits_{i\in I}(x^{\ast })^{i}\leq \tsum\limits_{i\in
I}e^{i}.$ $\ \ \ \ \ \ \ \ \ \ \ \ \ \ \ \ \ \ \ \ \ \ \ \ \ \ \ \ \ \ \ \ \
\ \ \ \ \ \ \ \ \ \ \ \ \square \medskip $

The second example is the abstract economy associated to an exchange economy
with two constraint set valued maps, the first one being the budget set
valued map and the second one being the consumption set that depends on
prices.\medskip

The third example follows the idea of an exchange economy which has, beyond
the budget set valued map, a second constraint set valued map $G_{i},$
defined by the delivery conditions as stated in the paper by
Correia-da-Silva and Herves-Beloso [6].

Let's assume that the set of the states of nature is $\Omega =\{1,2,...,m\}$%
, the future prices are $p_{1},p_{2},...,p_{m}\in IR_{+}^{l}$ and that each
agent $i$ has a signal $f_{i}:\Omega \rightarrow Y_{i}$ such that $%
f_{i}(s)=f_{i}(s^{\prime })$ if $s$ and $s^{\prime }$ are states that cannot
be distinguished. The agent $i$ chooses a portfolio $y(s)$ in the following
way:

$p_{s}y(s)\leq p_{s}y(s^{\prime })$ for all $s^{\prime }$ such that $%
f_{i}(s^{\prime })=f_{i}(s).$

The set valued map $G_{i}:X\times \Delta \rightrightarrows IR^{lm}$ is
defined by $G_{i}(x,p)=\{y\in IR^{lm}:p_{s}y(s)\leq p_{s}y(s^{\prime })$ for
all $s^{\prime }$ such that $f_{i}(s^{\prime })=f_{i}(s)\}.\medskip $

\noindent \textbf{4.3 The Existence of Equilibria in Locally Convex Spaces}

As an application of the fixed point Theorem 3.1, we have the following
result.\medskip

\noindent \textbf{Theorem 4.1} \textit{Let }$\Gamma =\left\{
X_{i},A_{i},B_{i},P_{i}\right\} _{i\in I}$\textit{\ be an abstract economy
such that for each }$i\in I,$\textit{\ the following conditions are
fulfilled:}

\textit{1)\ }$X_{i}$\textit{\ is a nonempty convex subset of a Hausdorff
locally convex topological vector space }$E_{i}\,$\textit{\ and }$D_{i}$%
\textit{\ is a nonempty compact convex subset of }$X_{i}$\textit{;}

\textit{2)\ for each }$x\in X:=\prod\limits_{i\in I}X_{i},$\textit{\ }$%
A_{i}\left( x\right) $\textit{\ and }$P_{i}(x)$\textit{\ are convex,\ }$%
B_{i}\left( x\right) $\textit{\ is nonempty, convex and }$A_{i}\left(
x\right) \cap P_{i}(x)\subset B_{i}(x)$\textit{;}

\textit{3) }$W_{i}:=\left\{ x\in X:A_{i}\left( x\right) \cap P_{i}\left(
x\right) \neq \emptyset \right\} $\textit{\ is open in }$X$\textit{.}

\textit{4)\ }$H_{i}:X\rightrightarrows X_{i}$\textit{\ defined by }$%
H_{i}\left( x\right) :=A_{i}(x)\cap P_{i}\left( x\right) $\textit{\ for each 
}$x\in X$\textit{\ is almost w-upper semicontinuous with respect to }$D_{i}$%
\textit{\ on }$W_{i}$\textit{\ and }$\overline{H_{i}^{V_{i}}}$\textit{\ is
convex nonempty valued for each open absolutely convex symmetric
neighborhood }$V_{i}$\textit{\ of }$0$\textit{\ in }$E_{i}$\textit{;}

5)\textit{\ }$B_{i}:X\rightrightarrows X_{i}$\textit{\ is almost w-upper
semicontinuous with respect to }$D_{i}$ and $\overline{B_{i}^{V_{i}}}$ 
\textit{is} \textit{convex nonempty valued for each open absolutely convex
symmetric neighborhood }$V_{i}$ \textit{of} $0$ \textit{in} $E_{i}$\textit{;}

6)\textit{\ for each }$x\in X$\textit{\ , }$x_{i}\notin \overline{\mathit{(}%
A_{i}\cap P_{i})}\left( x\right) $\textit{;}

\textit{Then there exists }$x^{\ast }\in D=$\textit{\ }$\prod\limits_{i\in
I}D_{i}$\textit{\ such that }$x_{i}^{\ast }\in \overline{B}_{i}\left(
x^{\ast }\right) $\textit{\ and }$(A_{i}\cap P_{i})(x^{\ast })=\emptyset $
for each $i\in I.\medskip $

\noindent \textit{Proof }Let $i\in I.$ \ By condition (3) we know that $%
W_{i} $ is open in $X.$

Let's define $T_{i}:X\rightrightarrows X_{i}$ by $T_{i}\left( x\right)
:=\left\{ 
\begin{array}{c}
A_{i}\left( x\right) \cap P_{i}\left( x\right) ,\text{ if }x\in W_{i}, \\ 
B_{i}\left( x\right) ,\text{ \ \ \ \ \ \ \ \ \ \ \ if }x\notin W_{i}%
\end{array}%
\right. $ for each $x\in X.$

Then $T_{i}:X\rightrightarrows X_{i}$ is a set valued map with nonempty
convex values. We shall prove that $T_{i}:X\rightrightarrows D_{i}$ is
almost w-upper semicontinuous with respect to $D_{i}$. Let \ss $_{i}$\ be a
basis of open absolutely convex symmetric neighborhoods of $0$ in $E_{i}$
and let \ss =$\tprod\limits_{i\in I}$\ss $_{i}.$

For each $V=(V_{i})_{i\in I}\in \tprod\limits_{i\in I}$\ss $_{i},$ for each $%
x\in X,$ let for each $i\in I$

$B^{V_{i}}(x):=(B_{i}\left( x\right) +V_{i})\cap D_{i}$,

$F^{V_{i}}(x):=((A_{i}\left( x\right) \cap P_{i}\left( x\right) )+V_{i})\cap
D_{i}$ and

$T_{i}^{V_{i}}(x)"=\left\{ 
\begin{array}{c}
F^{V_{i}}(x),\text{ if }x\in W_{i}, \\ 
B^{V_{i}}(x),\text{\ if }x\notin W_{i}.%
\end{array}%
\right. $

For each open set $V_{i}^{\prime }$ in $D_{i}$, the set

$\left\{ x\in X:\overline{T_{i}^{V_{i}}}\left( x\right) \subset
V_{i}^{\prime }\right\} =$

$=\left\{ x\in W_{i}:\overline{F^{V_{i}}}(x)\subset V_{i}^{\prime }\right\}
\cup \left\{ x\in X\smallsetminus W_{i}:\overline{B^{V_{i}}}(x)\subset
V_{i}^{^{\prime }}\right\} $

$=\left\{ x\in W_{i}:\overline{F^{V_{i}}}(x)\subset V_{i}^{^{\prime
}}\right\} \cup \left\{ x\in X:\overline{B^{V_{i}}}(x)\subset V_{i}^{\prime
}\right\} .$

According to condition (4), the set $\left\{ x\in W_{i}:\overline{F^{V_{i}}}%
(x)\subset V_{i}^{\prime }\right\} $ is open in $X$. The set $\left\{ x\in X:%
\overline{B^{V_{i}}}(x)\subset V_{i}^{\prime }\right\} $ is open in $X$
because $\overline{B^{V_{i}}}$ is upper semicontinuous.

Therefore, the set $\left\{ x\in X:\overline{T_{i}^{V_{i}}}\left( x\right)
\subset V_{i}^{\prime }\right\} $ is open in $X.$ It shows that $\overline{%
T_{i}^{V_{i}}}:X\rightrightarrows D_{i}$ is upper semicontinuous. According
to Theorem 3.1, there exists $x^{\ast }\in D=$ $\prod\limits_{i\in I}D_{i}$
such that $x^{\ast }\in \overline{T}_{i}\left( x^{\ast }\right) ,$ for each $%
i\in I.$ By condition (5) we have that $x_{i}^{\ast }\in \overline{B}%
_{i}\left( x^{\ast }\right) $ and $(A_{i}\cap P_{i})(x^{\ast })=\emptyset $
for each $i\in I.$ $\ \ \ \ \ \ \ \ \ \ \ \ \ \ \ \ \ \ \ \ \ \ \ \ \ \ \ \
\ \ \ \ \ \ \ \ \ \ \ \ \ \ \ \ \ \ \ \ \ \ \ \ \ \ \ \ \ \ \ \ \ \square
\medskip $

\noindent \textit{Example}\textbf{\ }\textit{4.1} Let $\Gamma =\left\{
X_{i},A_{i},B_{i},P_{i}\right\} _{i\in I}$\textit{\ }be an abstract economy,
where $I=\{1,2,...,n\},$ $X_{i}:=[0,4]$ be a compact convex choice set, $%
D_{i}:=[0,2]$ for each $i\in I$ and $X:=\prod\limits_{i\in I}X_{i}$.

Let the set valued maps $A_{i},B_{i},P_{i}:X\rightrightarrows X_{i}$ be
defined as follows:

for each $(x_{1},x_{2},...,x_{n})\in X,$

$A_{i}(x):=\left\{ 
\begin{array}{c}
\lbrack 1-x_{i},2]\text{ if }x\in (0,\frac{1}{2})^{n}; \\ 
\lbrack 1-x_{i},2)\text{ if }x\in \lbrack \frac{1}{2},1)^{n}; \\ 
\lbrack 3,4]\text{ \ \ \ \ \ \ \ \ \ \ \ if \ \ \ \ \ }x=0; \\ 
\lbrack 0,\frac{1}{2}]\text{, \ \ \ \ \ \ \ \ \ \ \ \ otherwise;}%
\end{array}%
\right. $

$P_{i}(x):=\left\{ 
\begin{array}{c}
\lbrack \frac{3}{2},2+x_{i}]\text{ if }x\in \lbrack 0,1)^{n}; \\ 
\lbrack 1,2]\text{, \ \ \ \ \ \ \ \ \ \ \ \ \ \ otherwise;}%
\end{array}%
\right. $

$B_{i}(x):=\left\{ 
\begin{array}{c}
\lbrack 0,2]\text{ if }x\in \lbrack 0,1); \\ 
\lbrack 3,4]\text{ \ \ \ if \ \ \ }x=0; \\ 
\lbrack 0,2)\text{,\ \ \ \ \ otherwise.}%
\end{array}%
\right. $

The set valued maps $A_{i},B_{i},P_{i}$ are not upper semicontinuous on $X.$

$A_{i}(x)\cap P_{i}(x):=\left\{ 
\begin{array}{c}
\lbrack \frac{3}{2},2]\text{ if }x\in (0,\frac{1}{2})^{n}; \\ 
\lbrack \frac{3}{2},2)\text{ if }x\in \lbrack \frac{1}{2},1)^{n}; \\ 
\phi \text{, \ \ \ \ \ \ \ \ \ \ \ \ \ \ otherwise.}%
\end{array}%
\right. $

$W_{i}:=\left\{ x\in X:A_{i}\left( x\right) \cap P_{i}\left( x\right) \neq
\emptyset \right\} =(0,1)^{n}$\textit{\ }is open in\textit{\ }$X.$

$\overline{\mathit{(}A_{i}\cap P_{i})}\left( x\right) :=\left\{ 
\begin{array}{c}
\lbrack \frac{3}{2},2]\text{ if }x\in \lbrack 0,1]^{n}; \\ 
\phi \text{, \ \ \ \ \ \ \ \ \ \ \ otherwise.}%
\end{array}%
\right. $

We notice that for each\textit{\ }$x\in X$\textit{\ , }$x_{i}\notin 
\overline{\mathit{(}A_{i}\cap P_{i})}\left( x\right) .$

We shall prove that $B_{i}$\ and $\mathit{(}A_{i}\cap P_{i})_{W_{i}}$ are
almost w-upper semicontinuous with respect to\textit{\ }$D_{i}=[0,2].$

On $W_{i},$

$\mathit{(}A_{i}\cap P_{i})\left( x\right) :=\left\{ 
\begin{array}{c}
\lbrack \frac{3}{2},2]\text{ if }x\in (0,\frac{1}{2})^{n}; \\ 
\lbrack \frac{3}{2},2)\text{ if }x\in \lbrack \frac{1}{2},1)^{n};%
\end{array}%
\right. ,$

$\mathit{(}A_{i}\cap P_{i})\left( x\right) +(-\varepsilon ,\varepsilon )=(%
\frac{3}{2}-\varepsilon ,2+\varepsilon )$ if $x\in (0,1)^{n};$

Let $\mathit{(}A_{i}\cap P_{i})^{V}\left( x\right) =(\mathit{(}A_{i}\cap
P_{i})\left( x\right) +(-\varepsilon ,\varepsilon ))\cap \lbrack 0,2],$
where $V=(-\varepsilon ,\varepsilon ).$

Then,

if $\varepsilon \in (0,\frac{3}{2}],$

$\mathit{(}A_{i}\cap P_{i})^{V}\left( x\right) =(\frac{3}{2}-\varepsilon ,2]$
if $x\in (0,1)^{n};$

if $\varepsilon >\frac{3}{2},$

$\mathit{(}A_{i}\cap P_{i})^{V}\left( x\right) =[0,2]$ if $x\in (0,1)^{n};$

Hence, for each $V=(-\varepsilon ,\varepsilon ),$ $\overline{\mathit{(}%
A_{i}\cap P_{i})^{V}}_{W_{i}}$ is upper semicontinuous and has nonempty
values.

$B_{i}\left( x\right) +(-\varepsilon ,\varepsilon )=\left\{ 
\begin{array}{c}
(-\varepsilon ,\text{ }2+\varepsilon )\text{ if },\text{ }x\in (0,1)^{n}; \\ 
(3-\varepsilon ,\text{ }4+\varepsilon )\text{\ \ \ if \ \ }x=0; \\ 
(-\varepsilon ,\text{ }2+\varepsilon )\text{ \ \ \ \ \ \ otherwise.}%
\end{array}%
\right. $

Let $B_{i}{}^{V}\left( x\right) =(B_{i}\left( x\right) +(-\varepsilon
,\varepsilon ))\cap \lbrack 0,2],$ where $V=(-\varepsilon ,\varepsilon ).$

Then,

if $\varepsilon \in (0,1],$

$B_{i}{}^{V}\left( x\right) =\left\{ 
\begin{array}{c}
\phi \text{ \ \ \ \ \ if \ }x=0; \\ 
\lbrack 0,2]\text{\ otherwise;}%
\end{array}%
\right. $

if $\varepsilon \in (1,3],$ $B_{i}{}^{V}\left( x\right) =\left\{ 
\begin{array}{c}
\lbrack 0,2]\text{ if }x\in \lbrack 0,1)^{n}; \\ 
(3-\varepsilon ,2]\text{ if }x=0; \\ 
\lbrack 0,2],\text{ \ \ \ \ \ \ otherwise.}%
\end{array}%
\right. $

and if $\varepsilon >3,$ $B_{i}{}^{V}\left( x\right) =[0,2]$ if $x\in X.$

Then, for each $V=(-\varepsilon ,\varepsilon ),$ $\overline{B_{i}^{V}}$ is
upper semicontinuous and has nonempty values.

Therefore, all hypotheses of Theorem 4.1 are satisfied, so that there exist
equilibrium points. For example, $x^{\ast }=\{\frac{3}{2},\frac{3}{2},...,%
\frac{3}{2}\}\in X$ verifies $x_{i}^{\ast }\in \overline{B}_{i}\left(
x^{\ast }\right) $ and $(A_{i}\cap P_{i})(x^{\ast })=\emptyset .$ \bigskip

Theorem 4.2 deals with abstract economies which have dual w-upper
semicontinuous pairs of set valued maps.

\noindent \textbf{Theorem 4.2 }\textit{Let }$\Gamma =\left\{
X_{i},A_{i},B_{i},P_{i}\right\} _{i\in I}$\textit{\ be an abstract economy
such that for each }$i\in I,$\textit{\ the following conditions are
fulfilled:}

1)\textit{\ }$X_{i}$\textit{\ is a nonempty convex subset of a Hausdorff
locally convex topological vector space }$E_{i}\,$and\textit{\ }$D_{i}$%
\textit{\ is a nonempty compact convex subset of }$X_{i}$\textit{;}

2)\textit{\ for each }$x\in X:=\prod\limits_{i\in I}X_{i},$\textit{\ }$%
P_{i}(x)\subset D_{i},$\textit{\ }$A_{i}\left( x\right) \cap P_{i}\left(
x\right) \subset B_{i}\left( x\right) $ \textit{and }$B_{i}\left( x\right) $%
\textit{\ is nonempty;}

3)\textit{\ the set }$W_{i}:=\left\{ x\in X:A_{i}\left( x\right) \cap
P_{i}\left( x\right) \neq \emptyset \right\} $\textit{\ is open in }$X$%
\textit{;}

4)\textit{\ the pair }$(A_{i\mid \text{cl}W_{i}},P_{i\mid \text{cl}W_{i}})$ 
\textit{is dual almost w-upper semicontinuous with respect to }$D_{i}$%
\textit{, }$B_{i}:X\rightrightarrows X_{i}$\textit{\ is almost w-upper
semicontinuous with respect to }$D_{i}$;

5) \textit{if} $T_{i,V_{i}}:X\rightrightarrows X_{i}$ \textit{is defined by} 
$T_{i,V_{i}}(x):=(A_{i}(x)+V_{i})\cap D_{i}\cap P_{i}(x)$ \textit{for each} $%
x\in X,$ \textit{then the} \textit{set valued maps} $\overline{B_{i}^{V_{i}}}
$ \textit{and} $\overline{T_{i,V_{i}}}$ \textit{are} \textit{nonempty} 
\textit{convex valued for each open absolutely convex symmetric neighborhood 
}$V_{i}$ \textit{of} $0$ \textit{in} $E_{i}$\textit{;}

6)\textit{\ for each }$x\in X$\textit{\ , }$x_{i}\notin \overline{P}%
_{i}\left( x\right) $\textit{;}

\textit{Then, there exists }$x^{\ast }\in D:=$\textit{\ }$\prod\limits_{i\in
I}D_{i}$\textit{\ such that }$x_{i}^{\ast }\in \overline{B}_{i}\left(
x^{\ast }\right) $\textit{\ and }$A_{i}\left( x^{\ast }\right) \cap
P_{i}\left( x^{\ast }\right) =\emptyset $\textit{\ for all }$i\in I.\medskip 
$

\noindent \textit{Proof} For each $i\in I,$ let \ss $_{i}$\ denote the
family of all open absolutely convex symmetric neighborhoods of zero in $%
E_{i}$ and let \ss $=\tprod\limits_{i\in I}$\ss $_{i}.$ For each $%
V=\tprod\limits_{i\in I}V_{i}\in \tprod\limits_{i\in I}$\ss $_{i},$ for each 
$i\in I,$ let

$B^{V_{i}}(x):=(B_{i}\left( x\right) +V_{i})\cap D_{i}$ for each $x\in X$ and

$S_{i}^{V_{i}}\left( x\right) :=\left\{ 
\begin{array}{c}
T_{i,V_{i}}(x),\text{ \ \ \ \ \ \ \ \ \ \ \ \ \ \ \ \ \ \ if }x\in W_{i}, \\ 
B_{i}^{V_{i}}(x),\text{ \ \ \ \ \ \ \ \ \ \ \ \ \ \ \ \ if }x\notin W_{i},%
\end{array}%
\right. $

$\overline{S_{i}^{V_{i}}}$ has closed values. Next, we shall prove that $%
\overline{S_{i}^{V_{i}}}:X\rightrightarrows D_{i}$ is upper semicontinuous.

For each open set $V^{\prime }$ in $D_{i}$, the set

$\left\{ x\in X:\overline{S_{i}^{V_{i}}}\left( x\right) \subset V^{\prime
}\right\} =$

$=\left\{ x\in W_{i}:\overline{T_{i,V_{i}}}(x)\subset V^{\prime }\right\}
\cup \left\{ x\in X\smallsetminus W_{i}:\overline{B_{i}^{V_{i}}}(x)\subset
V^{\prime }\right\} $

=$\left\{ x\in W_{i}:\overline{T_{i,V_{i}}}(x)\subset V^{\prime }\right\}
\cup \left\{ x\in X:\overline{B_{i}^{V_{i}}}(x)\subset V^{\prime }\right\} .$

We know that the set valued map $\overline{T_{i,V_{i}}}(x)_{\mid W_{i}}:$ $%
W_{i}\rightrightarrows D_{i}$ is upper semicontinuous. The set $\left\{ x\in
W_{i}:\overline{T_{i,V_{i}}}(x)\subset V^{\prime }\right\} $ is \ open in $%
X. $ Since $\overline{B_{i}^{V_{i}}}(x):X\rightrightarrows D_{i}$ is upper
semicontinuous, the set $\{x\in X:\overline{B_{i}^{V_{i}}}(x)\}\subset
V^{\prime }$ is open in $X$ and therefore, the set $\left\{ x\in X:\overline{%
S_{i}^{V_{i}}}\left( x\right) \subset V^{\prime }\right\} $ is open in $X$.
It proves that $\overline{S_{i}^{V_{i}}}:X\rightrightarrows D_{i}$ is upper
semicontinuous. According to Himmelberg's Theorem, applied for the set
valued maps $\overline{S_{i}^{V_{i}}},$ there exists a point $x_{V}^{\ast
}\in D=$ $\prod\limits_{i\in I}D_{i}$ such that $(x_{V}^{\ast })_{i}\in
S_{i}^{V_{i}}\left( x_{V}^{\ast }\right) $ for each $i\in I.$ By condition
(5), we have that $(x_{V}^{\ast })_{i}\notin \overline{P_{i}}\left(
x_{V}^{\ast }\right) ,$ hence, $(x_{V}^{\ast })_{i}\notin \overline{%
A_{i}^{V_{i}}}\left( x_{V}^{\ast }\right) \cap \overline{P_{i}}\left(
x_{V}^{\ast }\right) $. \newline
We also have that cl Gr $(T_{i,V_{i}})\subseteq $ cl Gr $(A_{i}^{V_{i}})\cap 
$cl Gr $P_{i}.$ Then $\overline{T_{i,V_{i}}}(x)\subseteq $ $\overline{%
A_{i}^{V_{i}}}(x)\cap \overline{P_{i}}\left( x\right) $ for each $x\in X.$
It follows that $(x_{V}^{\ast })_{i}\notin \overline{T_{i,V_{i}}}%
(x_{V}^{\ast }).$ Therefore, $(x_{V}^{\ast })_{i}\in \overline{B^{V_{i}}}%
\left( x_{V}^{\ast }\right) .$

For each $V=(V_{i})_{i\in I}\in \tprod\limits_{i\in I}$\ss $_{i},$ let's
define $Q_{V}=\cap _{i\in I}\{x\in D:x\in \overline{B^{V_{i}}}\left(
x\right) $ and $A_{i}\left( x\right) \cap P_{i}\left( x\right) =\emptyset
\}. $

$Q_{V}$ is nonempty since $x_{V}^{\ast }\in Q_{V},$ and it is a closed
subset of $D$ according to (3). Then, $Q_{V}$ is nonempty and compact.

Let \ss =$\tprod\limits_{i\in I}$\ss $_{i}.$ We prove that the family $%
\{Q_{V}:V\in \text{\ss }\}$ has the finite intersection property.

Let $\{V^{(1)},V^{(2)},...,V^{(n)}\}$ be any finite set of $\text{\ss\ }$and
let $V^{(k)}=\underset{i\in I}{\tprod }V_{i}^{(k)}{}_{i\in I}$, $k=1,...,n.$
For each $i\in I$, let $V_{i}=\underset{k=1}{\overset{n}{\cap }}V_{i}^{(k)}$%
, then $V_{i}\in \text{\ss }_{i};$ thus $V\in \underset{i\in I}{\tprod }%
\text{\ss }_{i}.$ Clearly $Q_{V}\subset \underset{k=1}{\overset{n}{\cap }}%
Q_{V^{(k)}}$ so that $\underset{k=1}{\overset{n}{\cap }}Q_{V^{(k)}}\neq
\emptyset .$

Since $D$ is compact and the family $\{Q_{V}:V\in \text{\ss }\}$ has the
finite intersection property, we have that $\cap \{Q_{V}:V\in \text{\ss }%
\}\neq \emptyset .$ Take any $x^{\ast }\in \cap \{Q_{V}:V\in $\ss $\},$ then
for each $V\in \text{\ss },$

$x^{\ast }\in \cap _{i\in I}\left\{ x^{\ast }\in D:x_{i}^{\ast }\in 
\overline{B^{V_{i}}}\left( x\right) \text{ and }A_{i}\left( x\right) \cap
P_{i}\left( x\right) =\emptyset )\right\} .$

Hence, $x_{i}^{\ast }\in \overline{B^{V_{i}}}\left( x^{\ast }\right) $ for
each $V\in $\ss\ and for each $i\in I.$ According to Lemma\emph{\ }2.2,\emph{%
\ }we have that\emph{\ } $x_{i}^{\ast }\in \overline{B_{i}}(x^{\ast })$ and $%
(A_{i}\cap P_{i})(x^{\ast })=\emptyset $ for each $i\in I.$ $\ \ \ \ \ \ \ \
\ \ \ \ \ \ \ \ \ \ \ \ \ \ \ \ \ \ \ \ \ \ \ \ \ \square $\medskip

We now introduce the following concept, which also generalizes the concept
of lower semicontinuous set valued maps.

\noindent \textbf{Definition 4.7} Let $X$ be a non-empty convex subset of a
topological linear space $E$, $Y$ be a non-empty set in a topological space
and $K\subseteq X\times Y.$

The set valued map $T:X\times Y\rightrightarrows X$ has the e-USCS-property
(e-upper semicontinuous selection property) on $K,$ if for each absolutely
convex neighborhood $V$ of zero in $E,$ there exists an upper semicontinuous
set valued map with convex values $S^{V}:X\times Y\rightrightarrows X$ such
that $S^{V}(x,y)\subset T(x,y)+V$ and $x\notin $cl $S^{V}(x,y)$ for every $%
(x,y)\in K$.$\medskip $

The following theorem is an equilibrium existence result for economies with
constraint set valued maps having e-USCS-property.

\noindent \textbf{Theorem 4.3} \textit{Let }$\Gamma
=(X_{i},A_{i},P_{i},B_{i})_{i\in I}$\textit{\ \ be an abstract economy,
where }$I$\textit{\ is a (possibly uncountable) set of agents such that for
each }$i\in I:$

(1)\textit{\ }$X_{i}$\textit{\ is a non-empty compact convex set in a
locally convex space }$E_{i}$\textit{;}

(2)\textit{\ }cl $B_{i}$\textit{\ is upper semicontinuous with non-empty
convex values;}

(3)\textit{\ the set }$W_{i}:$\textit{\ }$=\left\{ x\in X\text{ / }\left(
A_{i}\cap P_{i}\right) (x)\neq \emptyset \right\} $ \textit{is open;}

(3)\textit{\ }cl $(A_{i}\cap P_{i})$\textit{\ has \ the e-USCS-property on }$%
W_{i}$\textit{.}

\textit{Then there exists an equilibrium point }$x^{\ast }\in X$ \textit{\
for }$\Gamma $\textit{,}$\ i.e.$\textit{, for each }$i\in I$\textit{, }$%
x_{i}^{\ast }\in \overline{B}_{i}(x^{\ast })$\textit{\ and }$A_{i}(x^{\ast
})\cap P_{i}(x^{\ast })=\emptyset $\textit{.\medskip }

\noindent \textit{Proof} For each\textit{\ }$i\in I$, let \ss $_{i}$ denote
the family of all open convex neighborhoods of zero in $E_{i}.$ Let $%
V=(V_{i})_{i\in I}\in \tprod\limits_{i\in I}$\ss $_{i}.$ Since cl $%
(A_{i}\cap P_{i})$ has the e-USCS-property on $W_{i}$, it follows that there
exists an upper semicontinuous set valued map $F_{i}^{V_{i}}:X%
\rightrightarrows X_{i}$ such that $F_{i}^{V_{i}}(x)\subset $cl $(A_{i}\cap
P_{i})(x)+V_{i}$ and $x_{i}\notin $cl $F_{i}^{V_{i}}(x)$ for each $x\in
W_{i} $.

Define the set valued map $T_{i}^{V_{i}}:X\rightrightarrows X_{i}$, by

$T_{i}^{V_{i}}(x):=\left\{ 
\begin{array}{c}
\text{cl }\{F_{i}^{V_{i}}(x)\}\text{, \ \ \ \ \ \ \ \ \ \ \ \ \ \ \ \ \ \ \
\ \ \ \ \ if }x\in W_{i}\text{, } \\ 
\text{cl }(B_{i}(x)+V_{i})\cap X_{i}\text{, \ \ \ \ \ \ \ \ \ if }x\notin
W_{i}\text{;}%
\end{array}%
\right. $

$B_{i}^{V_{i}}:X\rightrightarrows X_{i},$ $B_{i}^{V_{i}}(x):=$cl $%
(B_{i}(x)+V_{i})\cap X_{i}=($cl $B_{i}(x)+$cl $V_{i})\cap X_{i}$ is upper
semicontinuous by Lemma 2.1\emph{.}

Let $U$ be an open subset of $\ X_{i}$, then

$U^{^{\prime }}:=\{x\in X$ $\mid T_{i}^{V_{i}}(x)\subset U\}$

\ \ \ =$\{x\in W_{i}$ $\mid T_{i}^{V_{i}}(x)\subset U\}\cup \{x\in
X\setminus W_{i}$ $\mid $ $T_{i}^{V_{i}}(x)\subset U\}$

\ \ \ =$\left\{ x\in W_{i}\text{ }\mid \text{cl }F_{i}^{V_{i}}(x)\subset
U\right\} \cup \left\{ x\in X\mid \text{ }(\text{cl }B_{i}(x)+\overline{V_{i}%
})\cap X_{i}\subset U\right\} $

\ \ \ 

$U^{^{\prime }}$ is an open set, because $W_{i}$ is open, $\left\{ x\in W_{i}%
\text{ }\mid \text{cl }F_{i}^{V_{i}}(x)\subset U\right\} $ is open since cl$%
F_{i}^{V_{i}}(x)$ is an upper semicontinuous map on $W_{i}.$ We have also
that the set $\left\{ x\in X\mid \text{ }(\text{cl }B_{i}(x)+\text{cl }%
V_{i})\cap X_{i}\subset U\right\} $ is open since $($cl $B_{i}+$cl $%
V_{i})\cap X_{i}$ is u.s.c. Then $T_{i}^{V_{i}}$ is upper semicontinuous on $%
X$ and has closed convex values.

Define $T^{V}:X\rightrightarrows X$ by $T^{V}(x):=\underset{i\in I}{\prod }%
T_{i}^{V_{i}}(x)$ for each $x\in X$.

$T^{V}$ is an upper semicontinuous set valued map and has also non-empty
convex closed values.

Since $X$ is a compact convex set, by Fan's fixed-point theorem [8], there
exists $x_{V}^{\ast }\in X$ such that $x_{V}^{\ast }\in T^{V}(x_{V}^{\ast })$%
, i.e., for each $i\in I$, $(x_{V}^{\ast })_{i}\in T_{i}^{V_{i}}(x_{V}^{\ast
})$. If $x_{V}^{\ast }\in W_{i},$ $(x_{V}^{\ast })_{i}\in $cl $%
F_{i}^{V_{i}}(x_{V}^{\ast })$, which is a contradiction.

Hence, $(x_{V}^{\ast })_{i}\in $cl $(B_{i}(x_{V}^{\ast })+V_{i})\cap X_{i}$
and $(A_{i}\cap P_{i})(x_{V}^{\ast })=\emptyset ,$ i.e. $x_{V}^{\ast }\in
Q_{V}$ where

$Q_{V}=\cap _{i\in I}\{x\in X:$ $x_{i}\in $cl $(B_{i}(x)+V_{i})\cap X_{i}$
and $(A_{i}\cap P_{i})(x)=\emptyset \}.$

Since $W_{i}$ is open, $Q_{V}$ is the intersection of non-empty closed sets,
therefore it is non-empty, closed in $X$.

We prove that the family $\{Q_{V}:V\in \underset{i\in I}{\tprod }\text{\ss }%
_{i}\}$ has the finite intersection property.

Let $\{V^{(1)},V^{(2)},...,V^{(n)}\}$ be any finite set of $\underset{i\in I}%
{\tprod \text{\ss }_{i}}$ and let $V^{(k)}=(V_{i}^{(k)})_{i\in I}$, $%
k=1,...n.$ For each $i\in I$, let $V_{i}=\underset{k=1}{\overset{n}{\cap }}%
V_{i}^{(k)}$, then $V_{i}\in \text{\ss }_{i};$ thus $V=(V_{i})_{i\in I}\in 
\underset{i\in I}{\tprod }\text{\ss }_{i}.$ Clearly $Q_{V}\subset \underset{%
k=1}{\overset{n}{\cap }}Q_{V^{(k)}}$ so that $\underset{k=1}{\overset{n}{%
\cap }}Q_{V^{(k)}}\neq \emptyset .$

Since $X$ is compact and the family $\{Q_{V}:V\in \underset{i\in I}{\tprod }%
\text{\ss }_{i}\}$ has the finite intersection property, we have that $\cap
\{Q_{V}:V\in \underset{i\in I}{\tprod }\text{\ss }_{i}\}\neq \emptyset .$
Take any $x^{\ast }\in \cap \{Q_{V}:V\in \underset{i\in I}{\tprod }\text{\ss 
}_{i}\},$ then for each $i\in I$ and each $V_{i}\in \text{\ss }_{i},$ $%
x_{i}^{\ast }\in $cl$(B_{i}(x^{\ast })+V_{i})\cap X_{i}$ and $(A_{i}\cap
P_{i})(x^{\ast })=\emptyset ;$ but then $x_{i}^{\ast }\in \overline{B}%
_{i}(x^{\ast })$ from Lemma 2.2\emph{\ }and $(A_{i}\cap P_{i})(x^{\ast
})=\emptyset $ for each $i\in I$ \ so that $x^{\ast }$ is an equilibrium
point of $\Gamma $ in X. $\ \ \ \ \ \ \ \ \ \ \ \ \ \ \ \ \ \ \ \ \ \ \ \ \
\ \ \ \ \ \ \ \ \ \ \ \ \ \ \ \ \ \ \ \ \ \ \ \ \ \ \ \ \ \ \ \ \ \ \ \ \ \
\ \ \square $\medskip

\noindent \textbf{5. Concluding Remarks}

We proved a fixed point theorem for the w-upper semicontinuous{\small \ }set
valued maps. We also obtained results concerning the existence of equilibria
for Yuan's model of an abstract economy without continuity assumptions.

The first author thanks the support by Research Grants ECO2012-38860-C02-02
(Ministerio de Economia y Competitividad), RGEA and 10PXIB300141PR (Xunta de
Galicia and FEDER).

\end{document}